\documentclass[12]{amsart}
\usepackage{amsmath, amssymb, verbatim, graphicx, amsthm}

\newtheorem{theorem}{Theorem}[section]

\newtheorem{proposition}{Proposition}[section]

\newtheorem{problem}{Problem}

\numberwithin{equation}{section}

\newcommand{\CH}{\mathcal H}

\newcommand{\CB}{\mathcal B}

\newcommand{\CU}{\mathcal U}
\newcommand{\R}{\mathbb R}
\newcommand{\Z}{\mathbb Z}
\newcommand{\T}{\mathbb T}

\newcommand{\diff}{{\rm Diff}}

\newcommand{\beq}{\begin{equation}}
\newcommand{\eeq}{\end{equation}}
\newcommand{\beqn}{\begin{equation*}}
\newcommand{\eeqn}{\end{equation*}}
\newcommand{\beqy}{\begin{eqnarray}}
\newcommand{\eeqy}{\end{eqnarray}}
\newcommand{\beqyn}{\begin{eqnarray*}}
\newcommand{\eeqyn}{\end{eqnarray*}}

\begin{document}

\title[SRB Entropy]{Behavior of the SRB Entropy Functional  in Families of Hyperbolic Attractors and Expanding Maps}

\author{  Miaohua Jiang  }

\address{
  Department of Mathematics, Wake Forest University, Winston Salem, NC 27109, USA}

\email{ jiangm@wfu.edu}


\begin{abstract}
We review the current state of the study of the SRB entropy functional over open components of families of uniformly hyperbolic systems, including Anosov systems, hyperbolic attractors,  expanding maps, and Markov transformations. We address problems concerning  the functional's  infimum,  supremum, and local extrema;  the gradient flow induced by the functional;  and the functional's behavior near the boundary of open components.
\end{abstract}

\maketitle

\large
\baselineskip 16 pt

\centerline{\bf\Large }

\section{Motivation, Definition, and Problem Formulation}

Motivated partly by the Gallavotti-Cohen Chaotic Hypothesis \cite{GC,G96,G06}, Ruelle initiated the study of the dependence of the SRB measure on the underlying  dynamical system,  a diffeomorphism on a closed Riemannian manifold possessing a uniformly hyperbolic attractor \cite{Ru97}. Ruelle proved that the SRB measure differentiably depends the system when the system is sufficiently smooth and derived the  first-order derivative formula of the SRB measure, i.e., the linear response function as well as the derivative formula of the entropy with respect to the SRB measure \cite{Ru03}, which we would call it the \emph{SRB entropy} or  the \emph{SRB entropy functional}.

This note contains  a collection of  properties of this functional for mainly transitive Anosov systems and expanding maps in several recent publications and a collection of problems for further exploration. Existing results on the behavior of the SRB entropy are sporadic and limited. We hope that this note can motivate more study of this topic.

\subsection{Definitions of the SRB measure and the SRB entropy functional}

We start with a slightly more general definition of the SRB measure for hyperbolic systems.  We assume that there is a $C^r, r > 1 $ diffeomorphism  $f_1$,   on a compact Riemannian manifold $M$, possessing an invariant  locally maximal hyperbolic set $\Lambda$:
 \[ \cap_{i \in \Z} f_1^{i} ( U)=\Lambda,\]
 where $U  \supset \Lambda$ is an open set. We also assume $f_1$ is transitive on $\Lambda$. Such a map is structurally stable: there exists an $\epsilon >0$ such that any other $C^{r}$ diffeomorphism $f$ of $M$ in $O^{C^r}_\epsilon (f_1)$, the $\epsilon$-neighborhood of $f_1$,  is topologically conjugate to $f_1$ via a bi-H\"older continuous homeomorphism $h_f$ on $\Lambda$:
\[ f \circ h_f = h_f \circ f_1.\]  Furthermore, the image $h_f(\Lambda)$, denoted by $\Lambda_f$,  is a locally  maximal hyperbolic set of  $f$.

Using a symbolic representation given by a Markov partition of $\Lambda_f$, for any H\"older continuous function $\varphi$ on $\Lambda_f$, we obtain both the \emph{ existence} and the \emph{uniqueness} of an invariant measure (an \emph{equilibrium state}) $\mu_\varphi$ on  $\Lambda_f$  that satisfies the variational principle:
\[ P(\varphi) = h_{\mu_\varphi}(f) + \int_{\Lambda_f} \varphi d \mu_\varphi=\sup_{\mu } \left\{  h_{\mu}(f) + \int_{\Lambda_f} \varphi d \mu \right\} ,\]
where $P(\varphi)$ is the topological pressure for the potential function $\varphi$ and $h_\mu(f)$ is the metric (measure-theoretic, or Kolmogorov-Sinai) entropy of $f$ with respect to $\mu$. For definitions of these terms as well as hyperbolic maps, we refer readers to standard reference books such as \cite{Bowen, Ru78, KH, Mane}.  For an overview of  earlier major results on the SRB measure  for various systems see \cite{Young}.

When we choose the potential function $\varphi$ to be $ - \log J^uf$, where $J^uf$ is the Jacobian of $f$ restricted to the unstable manifold of $f$, the corresponding equilibrium state  is called the \emph{Sinai-Ruelle-Bowen (SRB) measure}  of $f$ and is denoted by $\rho_f$. The entropy of $f$ with respect to the SRB measure $\rho_f$ will be  denoted by $\CH(f)$ in this note and called the \emph{SRB entropy} for brevity.  The topological entropy of $f$, $h_{\rm top}(f)$ is defined to be the topological pressure for the zero potential, $h_{\rm top}(f)=P(0)$, or, the entropy of $f$ with respect to the equilibrium state for the zero potential function, i.e.,  the measure of maximal entropy.

We see that the SRB entropy $\CH(f)$ is  a functional defined for all diffeomorphisms within the neighborhood $O^{C^r}_\epsilon(f_1)$.  Due to the continuity of the potential function $- \log J^u f$ on $f$ in $C^r$-topology, we see that the SRB functional
 \[\CH(f) = P(-\log J^uf) + \int \log J^u f  d \rho_f\] is   continuous  in $f$ in $C^r$-topology   since  the topological pressure is an analytic function over the space of H\"older continuous functions and the equilibrium state $\mu_\varphi$ is the first-order derivative of the topological pressure. The previous statement is a little overly simplified since we will also need the continuity of the conjugate map $h_f$ in $f$, which can be deduced from the continuity of stable and unstable manifolds of $f$ in $C^{r}, r >1$ topology.  The continuity and the differentiability of the SRB entropy functional $\CH(f)$ is proved via the conjugating map $h_f$:
   \[\CH(f) = P(-\log J^uf\circ h_f) + \int_\Lambda  \log J^u f\circ h_f\    d h^*_f(\rho_f).\]

Now we extend the domain of the SRB entropy functional $\CH(f)$ in the following way: Since every map in the neighborhood $O^{C^r}_\epsilon(f_1)$ is also structurally stable over its locally maximal hyperbolic set, we can extend the domain of $\CH(f)$ by taking union of all such open $\epsilon$-neighborhoods with nonempty intersections. We denote this union by $\CU^{C^r}(f_1)$, which is a path-connected open component of diffeomorphisms of $M$. Alternatively, we can define $\CU^{C^r}(f_1)$ in the following way:  $f \in \CU^{C^r}(f_1)$ if and only if there exist $\{\epsilon_i, i= 1,\cdots ,k\}, \{f_i\}, i =1,2 ,\cdots, k$ such that $   f_{i+1} \in O^{C^r}_{\epsilon_i}(f_{i}), i=k-1, k-2, \cdots, 1 $ and  every map in $       O^{C^r}_{\epsilon_i}(f_{i}), i=1,\cdots, k$ is topologically conjugate to $f_{i}$.

The value of the SRB entropy of a map  is invariant under a differentiable (smooth) conjugacy since the potential functions for smoothly conjugating maps belong to the same homologous class \cite{Bowen, Ru78}.  Thus, we can consider the SRB entropy as a functional defined on the space of smooth conjugating classes.

\subsection{Problems of Interest}

We now  formulate some problems concerning the SRB entropy functional over the domain $\CU^{C^r}(f_1)$.
Most of these problems are waiting to be explored and solutions to them are dependent on  individual families of systems. The same problems can be formulated for families of $C^r$ expanding maps or piecewise  expanding maps on Riemannian manifolds as long as maps in these families are structurally stable and the SRB measure exists uniquely.  Detailed description of known results and open problems will be presented in Sections \ref{sec.Anosov} and \ref{sec.expanding}.

\subsubsection{Range of the  functional}

 The  first set of problems  concerns  {  the range of values   of the SRB entropy functional}  $\CH(f), f \in \CU^{C^r}(f_1)$:

$\textbf{1.}$    {\sl  What are the values of $\inf_{f \in \CU^{C^r}(f_1)} \CH(f)$ and $\sup_{f \in \CU^{C^r}(f_1)} \CH(f)$} ?

Since $\rho_f$ is an equilibrium state of a transitive uniformly hyperbolic system, its entropy is always positive, so we have $\CH(f) \in (0, h_{top}(f)]$, where $h_{top}(f)$ is the topological entropy of $f$. Naturally, we ask whether $\inf_{f \in \CU^{C^r}(f_1)} \CH(f)=0$ and $\sup_{f \in \CU^{C^r}(f_1)} \CH(f) = h_{\rm top}(f)$.

The  path-connected open component $\CU^{C^r}(f_1)$ has an inherent Banach manifold structure. Each $C^r$ neighborhood of $f \in \CU^{C^r}(f_1)$ can be identified with an open $C^r$ neighborhood of the zero section of the Banach space consisting of $C^r$ vector fields on the manifold $M$.  Thus,
the SRB entropy $\CH(f)$  is a functional defined on an infinite dimensional Banach manifold. Existence of local extreme is a possibility:

$\textbf{2.}$  {\sl   Does the functional $\CH(f)$ have any local extreme in $\CU^{C^r}(f_1)$?}

 Since the space is of infinite dimension, the chance to have a local maximum smaller than the global maximum seems small since we can perturb the map in so many different ways. We will mention the potential physical meaning  of absence of a nontrivial local maximum when we discuss the next set of problems. 

\subsubsection{Gradient flow of the  functional}
  The  second set of questions  arises from the consideration of a {   gradient flow of the SRB entropy functional } in  $\CU^{C^r}(f_1)$.

When $r \ge 3$, the SRB entropy $\CH(f)$ is proven to be Fr\`echet differentiable in $f$: it is $C^{r-2}$ \cite{Ru97}. Let $n$ be the dimension of the manifold $M$.   The Sobolev embedding theorem \cite{AF, He, CJ25}  implies that  when $k > r + \frac{n}{2}$, the completion of the family of $C^k$ maps in the $C^r$ neighborhood of $f_1$ under the Sobolev norm of order $k$ (the $k$-th derivatives of $f$  belong to $L^2(M)$)  forms a
 Hilbert submanifold structure in $\CU^{C^r}(f_1)$.  We denote this Hilbert submanifold by $\CU^{H^k}(f_1)$. Within this  Hilbert submanifold, the gradient vector of  $\CH(f)$ is well-defined and the gradient vector field is Lipschitz continuous. Thus,
 the SRB entropy induces a gradient flow in $\CU^{H^k}(f_1) \subset \CU^{C^r}(f_1)$.

$\textbf{1.}$  {\sl Does the gradient flow exist globally?}

Since the Hilbert manifold $\CU^{H^k}(f_1)$ is not compact, the global existence of  the gradient flow is not guaranteed even thought the flow exists locally.

  If the answer is affirmative to the previous question, then we can further ask:   
  
    $\textbf{2.}$ {\sl  Do  trajectories of the SRB entropy gradient flow  converge as $t\to \pm \infty$?}

   In thermodynamics, the process for an isolated nonequilibrium system to evolve to an equilibrium is a diffusive process governed by the Boltzmann equation.  It inspires us to ask the next question:
   
   $\textbf{3.}$ {\sl   Is the process induced by the SRB entropy functional also a diffusive process?  Can the diffusive process induced by the SRB entropy be described by differential equations such as the ones in \cite{JKO,Maas}?}

\subsubsection{Boundary behavior of the  functional}

 The third set of questions  concerns  the value of the SRB entropy near the {\bf  boundary of the open component}  $\CU^{C^r}(f_1)$, (or on the boundary , if the SRB measure exists for a map on the boundary). These questions are motivated by recent work in \cite{HJJ1, HJJ2, YJ} as well as \cite{Luzz}.

Let's $\CB^r(f_1)$ denote the boundary set of $\CU^{C^r}(f_1)$ in the space of $C^r$ diffeomorphisms. A $C^r$ diffeomorphism $g \in \CB^r(f_1)$ if and only if it is not in  $\CU^{C^r}(f_1)$ but its $\epsilon$ neighborhood $O^{C^r}_\epsilon (g)$ intersects $\CU^{C^r}(f_1)$ for any $\epsilon >0$.  Note that when $r_2 > r_1 > 1$,  $\CU^{r_2}(f_1) \subset  \CU^{r_1}(f_1)$   and  $\CB^{r_2}(f_1) \subset  \CU^{r_1}(f_1)$.
We may also consider boundary sets in even broader sense: for example, a diffeomorphism $g$ with $L^1$ derivative  is in the boundary set $\CB^{H^1}(f_1)$ if it is not in $\CU^{C^r}(f_1)$ but any $H^1$ neighborhood of $g$ intersects $\CU^{C^r}(f_1)$.

  $\textbf{1.}$ {\sl Given $g \in \CB^{r} (f_1)$. Under what conditions does the SRB measure exist for $g$?  If $f_n \to g$ in $C^r$ topology, does the SRB entropy $\CH(f_n)$ converge to $\CH(g)$? }

    Let $\overline{\CH}_\epsilon(g) = \sup_{ d^{C^r}(f,g)\le \epsilon}   \CH(f)  $ and $\underline{\CH}_\epsilon(g) = \inf_{ d^{C^r}(f,g)\le \epsilon}   \CH(f)  $.

    $\textbf{2.}$ {\sl What is the difference $\overline{\CH}_\epsilon(g) -  \underline{\CH}_\epsilon(g) $ and how does this difference converge as $\epsilon \to 0$ for any $g \in \CB^r(f_1)$?}

The rest of the article is divided into two sections:  a section on hyperbolic diffeomorphisms  and another on expanding endomorphisms. We review the recent progresses concerning the behavior of the SRB entropy functional in   path-connected open components  and propose open problems for further exploration.

\section{ SRB entropy in a path-connected component of transitive uniformly hyperbolic attractors}\label{sec.Anosov}

This section is a collection of results on the possible values of the SRB entropy functional in families of transitive   hyperbolic attractors including Anosov maps,  published in  \cite{HJJ1,HJJ2, SVV, J24}

Theorem \ref{thm.attractor} concerns the infimum of the SRB entropy over $\CU^{C^r}(f_1)$ when $f_1$ possesses a uniformly hyperbolic attractor $\Lambda$:  there exists an open neighborhood $U$ of $\Lambda0$ such that
\[ \Lambda = \cap_{i = 0}^\infty f_1^{i}(U).\]
In this case, the topological pressure of the potential function $\varphi= -\log J^u  f$ is equal to zero for any $f \in \CU^{C^r}(f_1), $  $r > 1$ . Thus, we have the SRB entropy formula
\[      \CH(f) =   \int_{\Lambda_f}  \log J^u  f d \rho_f,\] where $\Lambda_f$ is the corresponding hyperbolic attractor of $f$ and $\rho_f$ is the SRB measure of $f$.

\begin{theorem}\label{thm.attractor}
\cite{HJJ1} Assume that $f_1  \in \diff^{C^r}(M)$ is a $C^r, r > 1$ diffeomorphism of a dimension $n\ge 2$ closed Riemannian manifold possessing a hyperbolic attractor $f $ and $\CU^{C^r}(f_1 )$ is the path-connected open component in
$\diff^{C^r}(M)$ containing $f_1$. Then there is $C^1$ path starting from $f_1$
\[ H=\{  f_t \in \CU^{C^r}(f_1) | 0<  t \le 1 \} \] such that   \[ \lim_{t \to 0^+}  \CH(f_t) = 0,\]
where $\CH(f_t)$ is the SRB entropy of $f_t$ on the corresponding hyperbolic attractor.
\end{theorem}

The idea of the proof is to perturb the map $f_1$ gradually near a fixed point $p$ (or a periodic point) so that the expansion rate goes to zero, i.e. $ \lim_{ t \to 0}   \log J^u  f_t  = 0$, at point $p$. We then estimate the integral
$ \int_{\Delta_{f_t}}  \log J^u  f_t d \rho_{f_t}$ over two regions: inside and outside of an small $\epsilon$-neighborhood of $p$. Outside this neighborhood, the integrant $\log J^u  f_t $ is bounded but the measure $\rho_{f_t}$'s density along the unstable manifolds is small while inside the $\epsilon$-neighborhood,  the integrand $\log J^u  f_t $ is nearly zero.
The limiting map $f_0$ is a diffeomorphism possessing an {\it almost} hyperbolic attractor (hyperbolic at every point other than the fixed point $p$) with a $\sigma-$finite SRB measure whose entropy is zero.

We observe that in the construction of the $C^1$-path $\{ f_t \}, t \in [0,1]$ in this theorem, the SRB measure $\rho_{f_t}$ changes and converges in weak$^*$ topology to a point mass measure as $t \to 0$.  What would it happen if we impose the condition that $\rho_{f_t}$  stays the same for all $t \in (0,1]$. Would this restriction become  a barrier to reduce  the SRB entropy  since we won't be able to break the entropy integral into two parts to have the desired estimates? In particular, what if we require that the diffeomorphism preserves a volume form of the Riemannian manifold $M$? The answer to this question was answered  in a subsequent paper by the same authors in the context of Anosov diffeomorphisms.

\begin{theorem}\label{thm.Anosov}
\cite{HJJ2} Assume that $f_1  \in \diff^{C^r}(M)$ is a $C^r, r > 1$ transitive Anosov diffeomorphism of a dimension $n\ge 2$ closed Riemannian manifold preserving a  smooth volume form and $\CU^{C^r}(f_1 )$ is the path-connected open component in the volume-preserving subspace of
$\diff^{C^r}(M)$ containing $f_1$. Then there is $C^1$ path in this component starting from $f_1$
\[ H=\{  f_t \in \CU^{C^r}(f_0) | 0<  t \le 1 \} \] such that   \[ \lim_{t \to 0^+}  \CH(f_t) = 0,\]
where $\CH(f_t)$ is the (SRB) entropy of $f_t$ with respect to the volume.
\end{theorem}

The discussion of the idea of the proof  of this theorem will be postponed after we state a similar theorem for the expanding maps in next section.

{\bf Remarks:} 

1.  The infimum of the SRB entropy problem remains a project incomplete. We do not know  whether the theorem holds in the cases of  (1) $f$ preserves a non-smooth SRB measure, i.e., an SRB measure without a continuous density function with respect to the Lebesgue measure.  (2) $f$ is an Axiom A diffeomorphism, where  the topological pressure is negative. In this case,  we have \[ \CH(f) = \int   \log J^u f d \rho_f + P( -\log \log J^u f) < \int   \log J^u f d \rho_f.\]
The perturbation described in the proof of Theorem \ref{thm.attractor} may still be able to reduce the value of
the integral  $\int   \log J^u f d \rho_f $ to nearly zero. Thus, the infimum of the SRB measure is still zero. But the details need to be verified since in the absence of stable foliations, the distortion estimates should be different.

2.  In the construction of the $C^1$-path $\{ f_t \}, t \in (0,1]$ in the proof of this theorem, the perturbation is not a small scale perturbation in $C^2$ topology. Indeed, as $ t \to 0$, the $C^2$ norm of the perturbation is not bounded.  If we impose the boundedness condition in $C^2$ norm of $f_t$,  we do not know whether the infimum is still zero.

3. On the other end of the spectrum of the value of the SRB entropy functional, the problem of determining the supremum of the SRB entropy has not been addressed.  We know that the entropy of any invariant measure does not exceed the topological entropy.  However, whether the supremum of the SRB entropy in $\CU^{C^r}(f_1)$  is the same as the topological entropy, in general, is unknown. For Anosov systems on a torus $\T^n$, if  $\CU^{C^r}(f_1)$  contains an Anosov automorphism, then the topological entropy is the maximum of the SRB entropy over $\CU^{C^r}(f_1)$.  We now know that in the family of Anosov maps on a high dimensional torus, there are path-connected open components that do not contain any automorphisms \cite{Gog, FG}.

For Anosov systems on a dimension two torus, every path-connected component  $\CU^{C^r}(f_0)$ contains an Anosov automorphism \cite{FG}. Also, the SRB entropy is proven not having any nontrivial local extrema. The following result is paraphrased  from \cite{SVV,J24}.

\begin{theorem}\label{thm.Anosov2}
 \cite{SVV, JL22}   Assume that $f_1  \in \diff^{C^r}(\T^2), r \ge 3 $ is a transitive Anosov diffeomorphism on a torus $\T^2$ preserving the Riemannian volume. Then, the entropy functional $\CH(f)$ over the open component $ \CU^{C^r}(f_1)$ does not have any  local minimum nor any nontrivial local  maximum. $\CH(f)$ reaches its maximum value only at diffeomorphisms that are smoothly conjugate to an automorphism.
\end{theorem}

The proof of this theorem relies on a detailed analysis of the consequences when the derivative operator of the map $ f   \to  \CH(f)$ is identically zero.  Similar results are unknown for Anosov systems on a torus of dimension higher than two.

\section{ SRB entropy functional in  families of  expanding maps}\label{sec.expanding}

We recall the definition of expanding maps on a closed Riemannian manifold $M$.  A $C^1$ map $f$ is called an expanding map on $M$ if there exist    constants $C>0$ and $\lambda > 1$ such that
the derivative operator of $f$, $Df$ satisfies the condition
\[ \| Df^n v \| \ge C \lambda^n \|v\|,   \forall v \in  T_x M, n \ge 1 .\]   The collection of all $C^r, r\ge 1$ expanding maps on $M$ is denoted by  ${\mathcal E}^r(M) $. We use ${\mathcal E}^r_\rho(M) $ to denote the subfamily of maps in ${\mathcal E}^r(M) $ that preserves a  volume $\rho$  on $M$.   Any $C^{1+\alpha}, \alpha > 0$ expanding map on $M$ preserves a volume $\rho$ which has a H\"older ($C^\alpha$)  continuous density function. Thus, when $f $ is at least $C^2$, we can always assume it is volume-preserving and its SRB measure is just the volume.

\subsection{Extreme values of the SRB entropy}

As we know that if a Riemannian manifold admits an expanding map, then its universal  cover is homeomorphic to $\R^n$ \cite{Grom}.
For expanding maps on a torus, it is homotopic to a unique automorphism within the family - an expanding map whose lift to the universal cover is a linear expanding map. Since the automorphism is Lebesgue measure preserving, so its SRB measure is also the measure of maximal entropy. Thus, the supremum of the entropy functional is the topological entropy which is attained at this automorphism.

In general, it is unknown to the author, whether every Riemannian manifold that admits an expanding map also admits an expanding map whose SRB measure and the measure of maximal entropy coincide.
Thus, in some results, we restrict  to  expanding maps on a torus or the circle.

 \begin{theorem} \label{expandinginf} \cite{HJJ2}
 For every  expanding map $f_1 \in {\mathcal E}^r_\rho(M), r > 2, $  there exists a $C^1$ path in ${\mathcal E}^r_\rho(M)$ starting from $f_1$:   $\{ f_t  \in  {\mathcal E}^r_\rho(M) | 0 < t \le 1 \} $  such that   $\lim_{t \to 0^+}  \CH(f_t)  = 0.$
 \end{theorem}

Proofs of Theorem \ref{thm.Anosov} and Theorem   \ref{expandinginf}    depend  on several ingredients: starting from a fixed point or a periodic point, the map is embedded into a hyperbolic flow. Then the flow  is slowed down using the technique of  Katok's perturbation \cite{Katok} for  hyperbolic flows to create a map that preserves a new smooth volume and has a smaller SRB entropy.  Applying the Dacorogna-Moser Lemma \cite{DM}, the map is again perturbed  with a smooth conjugacy to preserve the original volume and the small SRB entropy.


The limiting map  as $t\to 0$ along the path, if exists,  may fail to have any hyperbolicity.  See an example later in this section.

Similar to the case of dimension two  measure-preserving Anosov maps, the SRB functional $\CH(f)$ does not have nontrivial local extremum in the family of $C^{1+\alpha}$ expanding maps on the circle $S^1$.

\begin{theorem}\label{thm:em} \cite{J21}
For any given map $ f \in   {\mathcal E}_\rho^{{1+\alpha}}(S^1) $,  there exists $\varphi $, a  $C^{1+\alpha}$ map on $S^1$ such that for any $\epsilon$ sufficient small,
$ f + \epsilon \varphi \in  {\mathcal E}_\rho^{{1+\alpha}}(S^1)           $ and
 \begin{equation} \frac{d}{d \epsilon} H(f + \epsilon \varphi ) \big|_{\epsilon=0} \not= 0,\end{equation}
with the only exception when $\psi\circ f \circ \psi^{-1} $ is the degree $n$ linear expanding map on $S^1$ where $\psi(x)$ is a $C^{1+\alpha}$ diffeomorphism of $S^1$ satisfying $\psi'(x)=\rho(x)$.
\end{theorem}

This theorem implies that if $0< \CH(f) < h_{top}(f)$, $\CH(f)$ cannot be a local extremum.

\subsection{Differentiability of the SRB entropy and the Gradient Flow}

In \cite{Ru97},  the SRB entropy is proven to be a  Fr\`etch differentiable functional over a $C^3$-family of transitive uniformly hyperbolic systems. Via the transfer operator method, this differentiability can be proved for the SRB entropy over $C^2$-family of expanding maps \cite{Li03, Ba, Ba18}.   The Sobolev embedding theorem \cite{AF, He} implies that the SRB functional is also Fr\`echet differentiable over a subspace of $C^3$-family of transitive uniformly hyperbolic systems  which is endowed with a Hilbert manifold structure. Thus, the SRB entropy defines a gradient flow over this subspace  \cite{CJ25}.

While the local existence of the gradient flow follows directly from the fact that the derivative of the SRB entropy with respect to the map is Lipschitz continuous,  the global existence, i.e,  starting from an initial system with the SRB entropy within the interval  $(0, h_{top})$, whether the local orbit can be extended in both directions of time indefinitely is unknown except in the case when the manifold is $S^1$.  The gradient flow's existence for $t \in (-\epsilon, \infty)$ is only proved in the one dimensional measure-preserving expanding map case.  The flow is conveniently projected to the space of the inverse branches of expanding maps. The proof is based on the Riesz representation of the derivative operator on the  space of the inverse branches.

\begin{theorem}\label{thm.global}

The SRB entropy functional $f \to H(f) = \int_0^1 \ln f'(x) dx$, defined in $ {\mathcal E}_\rho^{{1+\alpha}}(S^1) $  induces a gradient flow on the space of derivatives of inverse map of $f$ under the $L^2$-norm.
This gradient flow exists globally for all $t \in [0, \infty)$ and every trajectory converges to the unique equilibrium which corresponds to the linear expanding map  at which $H(f)$ attains its maximum value.
 \end{theorem}

One natural question to ask is what happens if we let $t\to -\infty$? Would the flow trajectory be a path for the entropy to decrease towards zero? This question leads to the study of the behavior of the SRB entropy at the boundary of an open component of uniformly hyperbolic or expanding maps.

\subsection{SRB entropy near or in the boundary of uniformly expanding maps }

Here we give a couple of examples to illustrate possible behaviors of the SRB entropy for maps near or in the boundary set of the space of uniformly hyperbolic or expanding maps and a recent result for expanding maps that are defined through Blaschke products.

\textbf{ {Example 1. }} \cite{HJJ2} The identity map on a torus  is in the $L^1$-boundary of $C^\infty$ expanding maps: within any  $L^1$-neighborhood of the identity map, there are $C^\infty$ uniformly expanding maps with an arbitrarily small SRB entropy.

\begin{theorem}\label{example1}
There is a \( C^1 \)-path \( \{f_t\}_{0 < t < 1} \) of \( C^r \), \( r \geq 2 \), orientation-preserving circle expanding endomorphisms of degree \( d \geq 2 \) such that each \( f_t \) preserves the Lebesgue measure and the infimum of the metric entropy \( \CH(f_t) \to 0 \) as \( t \to 0^+ \). Moreover, \( f_t(x)  = \mathrm{id} + o(t) \) for all \( x \in S^1 \)  and $f_t \to  \mathrm{id}$ in $L^1(S^1)$.
\end{theorem}

The construction is based on smoothing  the following piecewise smooth map from
$[-\frac{1}{2}, \frac{1}{2}] \to \left[-\frac{1}{2}, \frac{1}{2}\right],$
\[L_t(x) =
\begin{cases}
\frac{1}{2t}\left(x + \frac{1}{2}\right), & x \in \left[-\frac{1}{2}, -\frac{1}{2}+ t\right], \\
\frac{x}{1-2t}, & x \in \left(-\frac{1}{2} + t, \frac{1}{2} - t\right), \\
\frac{1}{2t}\left(x - \frac{1}{2}\right), & x \in \left[\frac{1}{2} - t, \frac{1}{2}\right].
\end{cases}
\]

Here we identify the circle with $[-\frac{1}{2}, \frac{1}{2}]$.
This map  preserves the Lebesgue measure on the circle and its metric entropy is
\[ \CH(L_t) = -2t \log(2t) - (1 - 2t) \log(1 - 2t)  \]
and  $\lim_{t \to 0^+}  \CH(L_t) = 0 $.

\bigskip

On the other hand, it is also easy to construct maps in the boundary  set  with their  SRB entropy close to the maximal value, i.e. the topological entropy.

\textbf{ {Example 2.}}  There exist maps of nearly maximal SRB entropy  in the $H^1-$boundary of $\mathcal{E}^{C^r}(f)$.

In this example,  we identify $[-t,  1-t], t \in (0, \frac{1}{8})$ with the circle. We define
\[
M_t(x) =
\begin{cases}
\gamma_1(x), & x \in \left[-t , t \right], \\
2x,  (\text{mod} 1) & x \in \left(t, \frac{1}{2}- t  \right) \cup \left(   \frac{1}{2}+ t, 1 -t  \right)  \\
\gamma_2(x) (\text{mod} 1), & x \in \left[\frac{1}{2} - t, \frac{1}{2}+ t\right],
\end{cases}
\]
where $\gamma_1(x)$ is an odd function satisfying the following conditions:
(1)   $\gamma_1(x) $ is   chosen  such that $M_t(x) $ is $C^\infty$ over $[-t, t]$ including the endpoints,
(2)   $1< \gamma_1'(x) < 3 $   and  (3)  $\gamma_1'(0 ) = 1$, and
 $\gamma_2(x)$ is  chosen  based on how $\gamma_1(x)$ is defined such that $M_t(x)$ is $C^\infty $  on the circle except at the point $x=\frac{1}{2}$ and preserves the Lebesgue measure. We note that $\gamma'(\frac{1}{2})=\infty$

The entropy of $M_t(x)$ with respect to the Lebesgue measure is
$\CH(M_t)=\int_{-t}^{1-t}  \ln M'_t(x) d x  .$
We have
$\CH(M_t) \ge  2 \int_{ \frac{1}{2} + t }^{1-t}  \ln M'_t(x) d x   = 2 (\frac{1}{2}- 2t ) \ln 2= (1-4 t) \ln 2.$
We see that all maps $M_t$  are (nonuniformly) expanding  maps with two critical points: $x=0, \frac{1}{2} $ and preserves the Lebesgue measure and the derivative $M'_t(x)$ can be chosen to belong to $L^2[-t, 1-t].$

It is unknown to the author whether there exists a measure-preserving non-uniformly expanding map of degree 2 such that its entropy reaches the maximum value $\ln 2$.


\subsubsection{A family of expanding maps whose SRB entropy is completely understood}

A global picture of the SRB entropy  is quite complicated even for the family of expanding maps on the unit circle . In his recent work,  Yunping Jiang  analyzed expanding maps defined by Blaschke products in complex variables and obtain a complete picture in the degree 2 case:  a bell-shaped surface.

In the space of expanding Blaschke products
\[ B(z)= e^{i 2 \alpha\pi  } \prod^d_{n=1} \frac{ z - a_n}{1 - \bar a_n z}\ \  \text{where}\  |a_n| < 1 \ \text{and} \ 0\le \alpha < 1 ,\]
each Blaschke product defines an expanding map with degree $d$ when it has a fixed point in the open disk $D:=\{ z \in \mathbb C: |z| <1\}$.

\begin{theorem}\label{Blaschke} \cite{YJ}
The space of smooth conjugacy classes of expanding Blaschke products can be identified with the unit open disk $D$  when $d=2$. The SRB entropy of  expanding Blaschke products defines an analytic function in $D$ with level curves that are circles: $\{  c \in \mathbb C: |c| = r, 0\le r < 1\}$. The SRB entropy decreases to zero as $r \to 1$ and it is a concave downward function in $D$.
\end{theorem}

Most results on the SRB entropy for one dimensional maps on the circle can be extended to Markov transformations of the unit interval $[0,1]$.  We refer readers to \cite{JL22} for more details.

\newpage

\section*{Acknowledgments}
{ In preparing this note, I benefited from many discussions with Fan Yang.
The example of large entropy in the boundary set of uniform expanding maps was suggested to me by Yun Yang. I am also thankful for the long-time collaboration with Huyi Hu and Yunping Jiang on the study of the SRB entropy.}


\begin{thebibliography}{AB99}


\bibitem[AF]{AF} R.A. Adama and J.J.F. Fournier, Sobolev Spaces, 2nd ed. 2003, Elsevier Science Ltd.

\bibitem[Ba]{Ba} V. Baladi, Positive Transfer Operators and Decay of Correlations, World Scientific, Singapore, New Jersey, London Hong Kong, 2000

 \bibitem[Ba18]{Ba18} V.   Baladi, Viviane,
Dynamical zeta functions and dynamical determinants for hyperbolic maps.
A functional approach
Ergeb. Math. Grenzgeb. (3), 68 [Results in Mathematics and Related Areas. 3rd Series. A Series of Modern Surveys in Mathematics]
Springer, Cham, 2018. xv + 291 pp.
	
\bibitem[Bow]{Bowen}  R. Bowen,
Equilibrium states and the ergodic theory of Anosov diffeomorphisms.
Second revised edition. With a preface by David Ruelle. Edited by Jean-René Chazottes
Lecture Notes in Math., 470
Springer-Verlag, Berlin, 2008. viii+75 pp.

\bibitem[CJ25]{CJ25}  J. Chen and M. Jiang,  Lipschitz Continuity and Formulas of the Gradient Vector of the SRB
  Entropy Functional, arXiv:2509.18596, 1-27

\bibitem[DM]{DM} B. Dacorogna  and J. Moser,
On a partial differential equation involving the Jacobian determinant.
Ann. Inst. H. Poincar\'e Anal. Non Lin\'eaire 7 (1990), no. 1, 1-26.


\bibitem[FG]{FG} F. T. Farrell and A.  Gogolev,
The space of Anosov diffeomorphisms.
J. Lond. Math. Soc. (2) 89 (2014), no. 2, 383- 396


\bibitem[G96]{G96} G. Gallavotti,   Chaotic hypothesis: Onsager reciprocity and fluctuation-dissipation theorem.
J. Statist. Phys. 84 (1996), no. 5-6, 899-925

\bibitem[G06]{G06}G. Gallavotti,   Entropy, thermostats, and chaotic hypothesis. Chaos 16 (2006), no. 4, 043114, 6 pp

\bibitem[GC]{GC} G. Gallavotti and E.G.D.  Cohen,  Dynamical ensembles in stationary states. J. Statist. Phys. 80 (1995), no. 5- 6, 931-970.



\bibitem[Gog]{Gog}  A. Gogolev,
Diffeomorphisms H\"older conjugate to Anosov diffeomorphisms.
Ergodic Theory Dynam. Systems 30 (2010), no. 2, 441- 456.


\bibitem[Grom]{Grom}  M. Gromov,
Groups of polynomial growth and expanding maps.
Inst. Hautes \'Etudes Sci. Publ. Math. No. 53 (1981), 53 - 73

\bibitem[He]{He} E. Hebey,    Nonlinear Analysis on Manifolds: Sobolev spaces and inequalities, AMS,       1999

\bibitem[HJJ1]{HJJ1}Hu, Huyi; Jiang, Miaohua; Jiang, Yunping Infimum of the metric entropy of hyperbolic attractors with respect to the SRB measure. Discrete Contin. Dyn. Syst. 22 (2008), no. 1-2, 215-234.


\bibitem[HJJ2]{HJJ2} Hu, Huyi; Jiang, Miaohua; Jiang, Yunping Infimum of the metric entropy of volume preserving Anosov systems. Discrete Contin. Dyn. Syst. 37 (2017), no. 9, 4767-4783.





\bibitem[Luzz]{Luzz} D.  Coates, Douglas and S. Luzzatto,
Persistent non-statistical dynamics in one-dimensional maps.
Comm. Math. Phys. 405 (2024), no. 4, Paper No. 102, 34 pp



\bibitem[JKO]{JKO} R. Jordan, D. Kinderlehrer, and F. Otto, The variational formulation of the Fokker–Planck equation,
  SIAM J. Math.  Anal. 29 (1) (1998) 1–17.


\bibitem[J12]{J12} M. Jiang,    Differentiating potential functions of SRB measures on hyperbolic attractors, Ergodic Theory Dynam. Systems 32(4)  (2012) , 1350 - 1369

\bibitem[J21]{J21} M. Jiang,    Chaotic hypothesis and the second law of thermodynamics. Pure Appl. Funct. Anal. 6 (2021), no. 1, 205 - 219


\bibitem[JL22]{JL22} M. Jiang and M. Lopez, SRB entropy of Markov Transformations, J. Stat. Physics, 188 (2022) No.3 Paper No. 24, 18 pp

\bibitem[J24]{J24} M. Jiang, Gradient flow of the Sinai-Ruelle-Bowen entropy.
Comm. Math. Phys. 405 (2024), no. 5, Paper No. 118, 22 pp.



\bibitem[YJ]{YJ}  Y. Jiang, Global graph of metric entropy on expanding Blaschke products.
Discrete Contin. Dyn. Syst. 41 (2021), no. 3, 1469 - 1482.

\bibitem[Kat]{Katok}  A.  Katok,
Bernoulli diffeomorphisms on surfaces.
Ann. of Math. (2) 110 (1979), no. 3, 529 - 547

\bibitem[KH]{KH} A. Katok and B.  Hasselblatt,
Introduction to the modern theory of dynamical systems.
With a supplementary chapter by Katok and Leonardo Mendoza. Encyclopedia of Mathematics and its Applications, 54. Cambridge University Press, Cambridge, 1995

\bibitem[KL]{KL} G. Keller and C. Liverani, Stability of spectrum for transfer operators, Ann. Scuola Norm. Sup. Pisa Cl. Sci. 28, 141-152, 1999

\bibitem[Li03]{Li03} C. Liverani, Carlangelo,
Invariant measures and their properties. A functional analytic point of view. (English summary) Dynamical systems. Part II, 185–237.
Pubbl. Cent. Ric. Mat. Ennio Giorgi [Publications of the Ennio de Giorgi Mathematical Research Center]
Scuola Normale Superiore, Pisa, 2003

\bibitem[Maa]{Maas} J. Maas, Gradient flows of the entropy for finite Markov chains, J. Func. Analysis 261 (2011) 2250-2292


\bibitem[Man]{Mane} R. Ma\~n\'e,   Ergodic theory and differentiable dynamics. Translated from the Portuguese by Silvio Levy. Ergebnisse der Mathematik und ihrer Grenzgebiete (3) [Results in Mathematics and Related Areas (3)], 8. Springer-Verlag, Berlin, 1987.


\bibitem[Ru78]{Ru78} D. Ruelle,  Thermodynamic formalism.
The mathematical structures of equilibrium statistical mechanics. Second edition
Cambridge Math. Lib.
Cambridge University Press, Cambridge, 2004. xx+174 pp.


\bibitem[Ru89]{Ru89}  D. Ruelle,   The thermodynamic formalism for expanding maps. Comm.
Math. Phys., 125(2):239-262, 1989

\bibitem[Ru97]{Ru97}  D. Ruelle,   Differentiation of SRB states. Comm. Math. Phys. 187 (1997),
no. 1, 227-241.


\bibitem[Ru03]{Ru03}  D. Ruelle,  Correction and complements: ``Differentiation of SRB states'' [Comm. Math. Phys. 187 (1997), no. 1, 227–241; MR1463827].
Comm. Math. Phys. 234 (2003), no. 1, 185 -190


\bibitem[SS]{SS}  M. Shub and D.Sullivan. Expanding endomorphisms of the circle
revisited. Ergodic Theory Dynam. Systems, 5(2):285-289, 1985.

\bibitem[SVV]{SVV} R. Saghin, P. Valenzuela-Her\'iquez, and   C.H. V\'asquez,   Regularity with respect to the parameter of Lyapunov exponents for diffeomorphisms with dominated splitting.
Mem. Amer. Math. Soc. 300 (2024), no. 1506, v+78 pp

\bibitem[Yo]{Young} L-S. Young,  What are SRB measures, and which dynamical systems have them? Dedicated to David Ruelle and Yasha Sinai on the occasion of their 65th birthdays. J. Statist. Phys. 108 (2002), no. 5 -6, 733-754.



\end{thebibliography}
\end{document}